\documentclass[12pt,leqno]{article}
\usepackage{amsmath, amsthm, amsfonts, amssymb, color}
\usepackage{mathrsfs}
\usepackage{lmodern}
\usepackage[T1]{fontenc}
\usepackage[utf8]{inputenc}
\usepackage[hidelinks,bookmarks=false]{hyperref}
\hypersetup{pdftitle={Remarks on spectral gaps on Gaussian path space},pdfauthor={Bo Wu and Zhilai Wen}}
\newtheorem{thm}{Theorem}[section]

\newtheorem{counterexample}[thm]{Counterexample}
\newtheorem{remark}[thm]{Remark}
\newtheorem{defn}[thm]{Definition}
\newcommand{\scr}[1]{\mathscr{#1}}
\definecolor{wco}{rgb}{0.5,0.2,0.3}

\numberwithin{equation}{section}

\newcounter{assumption}

\title{{\bf Remarks on spectral gaps on Gaussian path space}\footnote{Supported in
		part by  NNSFC (12571154).} }

\author{{\bf Bo Wu and Zhilai Wen}\\\\
	\footnotesize { School  of Mathematical Sciences, Fudan
		University, Shanghai 200433, China}\\
{\footnotesize E-mail: wubo@fudan.edu.cn}}

\date{}

\begin{document}

\maketitle

	\newcommand{\R}{\mathbb R}
	\newcommand{\EE}{\mathbb E}
	\newcommand{\Z}{\mathbb Z}
	\renewcommand{\H}{\mathbb H}
	\newcommand{\HH}{\scr H}
	\newcommand{\<}{\langle}
	\renewcommand{\>}{\rangle}
	\newcommand{\pp}{\partial}
	\renewcommand{\d}{\text{\rm{d}}}
	\newcommand{\loc}{\text{\rm{loc}}}
	\newcommand{\D}{\scr D}
	\newcommand{\E}{\scr E}
	\newcommand{\ess}{\text{\rm{ess}}}
	\newcommand{\beg}{\begin}
	\newcommand{\F}{\scr F}
	\newcommand{\Ric}{\text{\rm{Ric}}}
	\newcommand{\Var}{\text{\rm{Var}}}
	\newcommand{\Ent}{\text{\rm{Ent}}}
	\newcommand{\Hess}{\text{\rm{Hess}}}
	\newcommand{\B}{\scr B}
	\newcommand{\e}{\text{\rm{e}}}
	\renewcommand{\b}{\mathbf b}
	\newcommand{\cut}{\text{\rm{cut}}}
	\renewcommand{\P}{\mathbb P}
	\newcommand{\C}{\scr C}
	\newcommand{\M}{\scr M}
	\newcommand{\X}{\scr X}
	\newcommand{\LL}{\scr L}
	\newcommand{\T}{\scr T}
	\newcommand{\A}{\mathbf A}
	\newcommand{\gap}{\mathbf{gap}}
	\renewcommand{\div}{\text{\rm div}}
	\newcommand{\Lip}{\text{\rm Lip}}
	\newcommand{\dist}{\text{\rm dist}}
	\newcommand{\supp}{\text{\rm supp}}
	\newcommand{\Cov}{\text{\rm Cov}}
	\newcommand{\Dom}{\text{\rm Dom}}
	\renewcommand{\Cap}{\text{\rm Cap}}
	\newcommand{\II}{{\mathbb I}}
	\newcommand{\sect}{\text{\rm sect}}
	\newcommand{\N}{\mathbb N}
	\newcommand{\K}{\scr K}
	\newcommand{\RR}{\scr R}
	\newcommand{\diag}{\operatorname{diag}}

\begin{abstract}
In this paper, we discuss the relationship between spectral gaps of Ornstein--Uhlenbeck/$L^2$ operators on Gaussian path space and the underlying Volterra kernel $F$. This framework includes fractional Brownian motion $B^H$ with Hurst index $H\in(0,1)$. For the fractional Wiener measure $\mu_H$, the logarithmic Sobolev inequality with respect to the Ornstein--Uhlenbeck Dirichlet form holds for $H\geq\frac12$, whereas the corresponding weak Poincar\'e inequality fails for $0<H<\frac12$. The damped/$L^2$ Dirichlet forms satisfy logarithmic Sobolev inequalities for every $H\in(0,1)$.

\end{abstract}

\vskip 1cm

\section{Introduction}

Let $(\Omega,\F, \mathbb P)$ be a complete probability space and let $\{B_t: t\geq 0\}$ denote a standard Brownian motion on $\mathbb R $. In this article, we mainly consider Volterra processes
$$ B_t^F := \int_0^t F(t,s)\d B_s, \quad t \geq 0,
$$
where $F:\mathbb R \times \mathbb R \to \mathbb R$ is a deterministic Volterra kernel.
We assume always that $B^F$ is a Gaussian process admitting a continuous version, which we use throughout the paper.
By its Volterra representation, $B^F$ is centered and has covariance
\begin{equation}
	\label{eq1.2}
	\mathbb E(B_t^F B_s^F)={\int_0^{s\wedge t}} F(t,r) F(s,r) \d r<\infty, \quad s,t\geq0.
\end{equation}
We refer to $B^F$ as an $F$-Gaussian process, or Volterra-Gaussian process. Examples include Brownian motion (corresponding to $F(t,s)=1$ for $s\leq t$) and fractional Brownian motion (corresponding to $F=K_H$), where, for $0<s<t$, the kernel $K_H$ is defined by
\begin{equation}
	\label{eqKF}
	K_H(t,s) =
	\begin{cases}
		c_H s^{\frac12 - H} \int_{s}^{t} (u-s)^{H-\frac32} u^{H-\frac12} \mathrm{d}u,
		& H \in \bigl(\frac12,1\bigr), \\
		\mathbf{1}_{[0,t]}(s),
		& H = \frac12, \\
		b_H \left( \frac{t(t-s)}{s} \right)^{H-\frac12}
		- b_H\bigl(H-\frac12\bigr) s^{\frac12 - H} \int_{s}^{t} (u-s)^{H-\frac12} u^{H-\frac32} \mathrm{d}u,
		& H \in \bigl(0,\frac12\bigr),
	\end{cases}
\end{equation}
where
$$\aligned &b_H = \sqrt{2H/((1-2H)\beta(1-2H,H+1/2))},\\
&c_H = \sqrt{ H(2H - 1) / \left( \beta(2 - 2H, H - 1/2) \right) }.\endaligned$$
Here $\beta(\cdot,\cdot)$ is the beta function (see \cite{BHOZ,N}).   

For any fixed $T>0$, let $\mu^F$ be the distribution of $B^F$ in the time interval $t\in [0,T]$.
 Then $\mu^F$ is a probability measure on the path space $W_T:= C([0,T];\R)$. We call $\mu^F$ the $F$-Gaussian Wiener measure. The class of cylinder functions $\F C_b^{L^2}$ on $W_T$ is given by
\begin{equation*}
	\begin{aligned}
	\F C_b^{L^2}=\bigg\{&G(\gamma)=g\bigg(\int_0^Tg_1(s,\gamma_s)\d s,\dots,\int_0^Tg_m(s,\gamma_s)\d s\bigg):~\gamma\in W_T, \\
	&g\in C_b^1(\R^{m}),g_i\in C^{0,1}([0,T]\times\R;\R),\ \|g'_{i,2}\|_\infty<\infty,1\le i\le m\bigg\},
	\end{aligned}
\end{equation*}
Here $C^{0,1}$ means that the function and its derivative with respect to the second variable are jointly continuous.

Let $\mathbb{H}$ be the standard Cameron--Martin space for $C([0,T];\mathbb{R})$, i.e.
\begin{equation*}
\begin{split}
&\mathbb{H}:=\Big\{h\in C([0,T]; \mathbb{R})\Big|
h~ \text{is absolutely continuous}\\
&~~~~~~~~~~~~~~~~~~~~~~~~~~~~~~~~h(0)=0, \|h\|^2_{\mathbb{H}}:=\int_0^T| h'(s)|^2\d s<\infty\Big\},
\end{split}
\end{equation*}
where $h'(s)$ is the derivative with respect to the time variable $s$. In
fact, $\mathbb{H}$ is
a separable Hilbert space equipped with the inner product $\<h,g\>_{\mathbb{H}}:=\int_0^T  h'(s) g'(s)\d s,\
h,g\in \mathbb H.$
For any $F\in\F C_b^{L^2}$ of the form
$F(\gamma):=g\bigg(\int_0^Tg_1(s,\gamma_s)\d s,\dots,\int_0^Tg_m(s,\gamma_s)\d s\bigg)$ and any
$h\in L^2([0,T];\d x)$, we define the directional derivative $D_h F$ as follows:
\begin{equation}\label{eqc1.1}
D_hF(\gamma):=\sum_{i=1}^{m}\partial_ig\int_0^T g'_{i,2}(s,\gamma_s)h(s)\d s,
\end{equation}
where $g'_{i,2}$ is the derivative with respect to the second variable.
By the Riesz representation theorem, for the above $F$ and $h$, there exist a Malliavin gradient $\nabla F\in \mathbb{H}$ and an $L^2$-gradient $D F$
such that
\[
\begin{aligned}
	\<\nabla F(\gamma),h\>_{\mathbb{H}} &= D_hF(\gamma), \quad h\in\mathbb{H}, \ \gamma \in W_T,\\
	\<D F(\gamma),h\>_{L^2([0,T];\d x)} &= D_hF(\gamma), \quad h\in L^2([0,T];\d x), \ \gamma \in W_T.
\end{aligned}
\]
In fact, it is straightforward to verify that the Malliavin gradient and the $L^2$-gradient are given, respectively, by
\[
\begin{aligned}
	(\nabla F(\gamma))'(r) &= \sum_{i=1}^m \partial_i g \int_r^T g'_{i,2}(s,\gamma_s)\d s,\quad r\in[0,T],\\
	(DF(\gamma))(r) &= \sum_{i=1}^m \partial_i g\, g'_{i,2}(r,\gamma_r),\quad r\in[0,T].
\end{aligned}
\]
To construct diffusion processes on Gaussian path space, firstly, we need to define two quadratic forms:
\begin{equation}\label{eqc1.2}
\begin{aligned}
\E_{OU}(F,G)&:=\int_{W_T}\<\nabla F,\nabla G\>_{\mathbb{H}}\d \mu^F,
&&F,G\in\F C_b^{L^2},\\
\E_{L^2}(F,G)&:=\int_{W_T}\<DF,DG\>_{L^2([0,T];\d x)}\d \mu^F,
&&F,G\in\F C_b^{L^2}.
\end{aligned}
\end{equation}
Under suitable assumptions, it is straightforward to show that the above quadratic forms are closable and their closures are quasi-regular Dirichlet forms(detailed refer to \cite{CSW}). The generator $\LL_{OU}$
 associated to the Dirichlet form $\E_{OU}$ and the generator $\LL_{L^2}$ associated to the Dirichlet form $\E_{L^2}$ are called the Ornstein-Uhlenbeck operator and $L^2$ operator respectively. For the $F$-Gaussian processes with non-constant $F$, Fan \cite{F,F15} proved the spectral gaps of $L^2$-operator for fractional Brownian motion. Chen--Sun--Wu \cite{CSW} established Log-Sobolev inequalities for Ornstein-Uhlenbeck/$L^2$-Dirichlet forms on certain reasonable conditions of density function $F$, but they do not cover the case of fractional Brownian motion with Hurst index $H<1/2$.

When $F\equiv1$ and $\R$ is replaced by a Riemannian manifold with bounded Ricci curvature, a spectral gap for the corresponding OU form was proved in \cite{Fang94}. Logarithmic Sobolev inequalities have been established by Aida--Elworthy \cite{AE95}, Capitaine--Hsu--Ledoux \cite{CHL97} and Hsu \cite{Hsu97}.
For some non-compact Riemannian manifolds with unbounded curvature, Wang \cite{Wang04} established weak Poincar\'e inequalities for manifolds with unbounded Ricci curvature, while Chen--Wu \cite{ChenWu14} derived Poincar\'e inequalities. For the $L^2$-Dirichlet form $\E_{L^2}$, Gourcy--Wu \cite{GW} established the logarithmic Sobolev inequality, which was later generalized to  broader settings in \cite{RWZZ, CWZZ}.

The paper is organized as follows. In Section 2, we introduce the stochastic integral with respect to $B^F$, together with the associated operators $K_F$ and $K_F^*$. Section 3 is devoted to the existence of spectral gaps for $\LL_{OU}$ and $\LL_{L^2}$. Moreover, under appropriate assumptions, we prove that the existence of a spectral gap of the Ornstein--Uhlenbeck operator is equivalent to the boundedness of $K_{F}^{*}$. Using this equivalence, we further show that the weak Poincar\'{e} inequalities for Ornstein-Uhlenbeck Dirichlet form with respect to fractional Wiener measure fail when the Hurst index satisfies $0<H<\frac{1}{2}$. Finally, we prove the Bismut--Elworthy--Li formula.

\section{Integration with respect to $F$-Gaussian processes}

This section introduces integration with respect to $B^F$ and the operators $K_F$ and $K_F^*$, which will be used in the quasi-invariance formula and the functional inequalities.
Throughout the paper, we impose the following standing assumptions on the Volterra kernel $F$:
\begin{enumerate}
\item[(S1)] $F$ is a measurable Volterra kernel on $[0,T]^2$, with $F(t,s)=0$ for $s>t$, and
\[
\phi_t:=F(t,\cdot)\in L^2([0,T])\quad\text{for every }t\in[0,T].
\]
\item[(S2)] The centered Gaussian process $B_t^F=\int_0^t F(t,s)\d B_s$ admits a continuous version; a sufficient condition is the following estimate for some constants $C>0$ and $\rho>0$:
\[
\EE\big[(B_t^F-B_s^F)^2\big]=\|\phi_t-\phi_s\|_2^2\leq C|t-s|^\rho\qquad\text{for all }s,t\in[0,T].
\]
\item[(S3)] For a.e. $s\in[0,T)$, the map $t\mapsto F(t,s)$ is locally absolutely continuous on $(s,T]$.

\item[(S4)] $
C_T:=\sup_{0\leq t\leq T}\|\phi_t\|_2^2<\infty.
$
\end{enumerate}
For centered Gaussian processes, continuity in probability implies mean-square continuity. Consequently, (S2) implies that $t\mapsto\phi_t$ is continuous in $L^2$, and (S4) follows by compactness. We retain (S4) to name the constant used below. The sufficient estimate in (S2) yields a continuous version by Gaussian moment estimates and Kolmogorov's criterion.

\subsection{The operators $K_F$ and $K_F^*$}

For the kernel $F$, define an operator $K_F$ on $L^2([0,T];\d x)$ by
$$ K_F h (t)  := \int_0^t F(t,s) h(s) \d s ,\quad h \in L^2([0,T];\d x),\quad t\in [0,T].
$$
Let $\mathbb H^F=K_F(L^2([0,T]))$, equipped with the quotient norm
\[
\|h^F\|_{\H^F}:=\inf\big\{\|u\|_2:K_Fu=h^F\big\}.
\]
Write $M=(\ker K_F)^\perp$. Each $h^F$ has a unique minimal preimage $u\in M$, and $K_F:M\to\H^F$ is an isometry. Also
\[
\|K_Fu\|_\infty\leq\sqrt{C_T}\|u\|_2,
\quad |K_Fu(t)-K_Fu(s)|\leq\|\phi_t-\phi_s\|_2\|u\|_2.
\]
Thus $\H^F$ is continuously embedded in $W_T$.

Define the transfer operator $K_F^*$ by the following increment duality:
We denote by $\mathcal S$ the vector space of step functions on $[0,T]$. For any $g\in \mathcal S$ and $h \in L^2([0,T];\d x)$,
$$ \int_0^T K_F^* g (t) h(t) \d t  = \int_0^T g(t)\d K_F h(t).
$$
If $g=\sum_{i=1}^n a_i1_{(s_{i-1},s_i]}$, the increment duality gives directly
\[
K_F^*g=\sum_{i=1}^n a_i(\phi_{s_i}-\phi_{s_{i-1}}).
\]
In particular, $K_F^*1_{[0,t]}=\phi_t$. By the standing local absolute-continuity assumption, the same formula yields
\begin{equation}\label{eqKF*}
K_F^*g(s)=g(s)F(T,s)+\int_s^T(g(t)-g(s))\partial_tF(t,s)\d t,
\quad g\in\mathcal S,
\end{equation}
where step functions are taken right-continuous away from their finitely many endpoints. For almost every $s$, the difference $g(t)-g(s)$ vanishes for $t$ sufficiently close to $s$, so local absolute continuity suffices for this representation.

The deterministic integrand space $\HH$ is the completion of $\mathcal S$ under the covariance seminorm, after quotienting by its null space. More precisely, set
$$\langle 1_{[0,a]},1_{[0,b]} \rangle_\HH := \int_0^T K_F^*1_{[0,a]}(t) K_F^*1_{[0,b]}(t) \d t = \int_0^{T} F(a,t)F(b,t)\d t$$
on step functions $\mathcal S$. For any $g\in\HH$, choose a sequence $(g_n)$ in $\mathcal S$ converging to $g$ in $\HH$ and define
$$K_F^*g := \lim_{n\to \infty} K_F^*g_n.
$$
The extension $K_F^*:\HH\to M$ is an isometric isomorphism. 

In contrast, let $\K_F^*$ denote the usual $L^2$ adjoint of $K_F$:
$$ \int_0^T \K_F^* g (t) h(t) \d t = \int_0^T g(t) K_Fh(t) \d t.
$$
where $g,h \in L^2([0,T];\d x)$. Fubini's theorem gives
$$\K_F^* g(t) = \int_t^T F(r,t)g(r) \d r.
$$

\begin{remark}
	Under (S4), $\|F\|_{L^2([0,T]^2)}^2\leq TC_T$. In particular,

	$(1)$ $K_F: L^2([0,T];\d x) \to L^2([0,T];\d x)$ is a bounded linear operator;
	
	$(2)$ $\K_F^*: L^2([0,T];\d x) \to L^2([0,T];\d x)$ is a bounded linear operator;
	
\end{remark}
Although $K_F^*:\HH\to M$ is an isometry, its restriction to step functions need not satisfy an estimate $\|K_F^*g\|_2\leq C\|g\|_2$. Boundedness of $K_F^*$ in the $L^2$ norm will always refer to this estimate. In particular, for the fractional Brownian kernel $F=K_H$, $K^*_F$ is a bounded operator for the Hurst index $H\geq \frac{1}{2}$, while $K^*_F$ does not satisfy this boundedness estimate for the Hurst index $H<\frac{1}{2}$. This distinction is central to the arguments concerning functional inequalities in Section 3.

\begin{remark}
	If $K_F$ is injective on $L^2([0,T];\d x)$, then $K_F^*(\mathcal S)$ is dense in $L^2([0,T];\d x)$. Indeed, for every $a\in[0,T]$, we have
	$$K_F^*1_{[0,a]}=F(a,\cdot).$$
	If $u\in L^2([0,T];\d x)$ is orthogonal to $K_F^*(\mathcal S)$, then
	$$0=\langle F(a,\cdot),u\rangle_{L^2([0,T];\d x)}=\int_0^a F(a,s)u(s)\,\d s=(K_Fu)(a),\quad a\in[0,T].$$
	The injectivity of $K_F$ therefore implies $u=0$, so that $\overline{K_F^*(\mathcal S)}$ is dense in $L^2([0,T];\d x)$.
\end{remark}

\begin{remark}
The representation \eqref{eqKF*} is used in \cite{AMN} for singular Gaussian kernels. For fractional Brownian motion, the relation between $\HH$ and spaces of functions or distributions is discussed in \cite{DU99,N03,PT,PT00}. 
\end{remark}

\subsection{Integration with respect to $F$-Gaussian processes}
We define deterministic Wiener integrals following Chapter 2 of \cite{BHOZ}.

\begin{defn}
	\label{defInt}
	The Wiener integral with respect to $B^F$ is the continuous linear extension to $\HH$ of the isometry initially specified by
	$$
	\begin{aligned}
		J^F: \HH &\to L^2(\P)
		\\ 1_{[0,t]} &\to B^F_t.
	\end{aligned}
	$$
\end{defn}
By Definition~\ref{defInt}, it follows that the integral with respect to step functions $\mathcal S$ is given by
$$J^F(\sum_{i=1}^n a_i1_{(t_i,t_{i+1}]})=\sum_{i=1}^n a_i(B^F_{t_{i+1}}-B^F_{t_i}).
$$
Now consider a general $ \psi\in \HH $. For any $\psi\in \HH $, there exists a sequence $ (\psi_n)_{n \in \mathbb{N}} \subset \HH $ such that every $ \psi_n\in \mathcal S $ and $\psi_n\to\psi$ in $\HH$. Hence we can define the integral of $\psi\in \HH $ with respect to $ B^F $ as
\[
\int_0^T \psi(s)\d B_s^F:=J^F(\psi) = \lim_{n \to \infty} J^F(\psi_n),
\]
where the limit is taken in $ L^2(\mathbb{P}) $.

The following transfer identity relates the Wiener integrals with respect to $B^F$ and $B$.
\begin{thm}
	\label{thmInt}
	For a fixed $T>0$, we have
	$$ \int_0^T \psi(s) \d B_s^F = \int_0^T K_F^* \psi(s) \d B_s,\quad \forall \psi \in \HH.
	$$
\end{thm}
\begin{proof}
	By the definition of $B_t^F$, we have
	$$ B_t^F = \int_0^T F(t,s) 1_{[0,t]}(s) \d B_s = \int_0^T K_F^* 1_{[0,t]}(s) \d B_s.
	$$
	and for $0\leq a<b\leq T$
	\begin{equation}
		\label{eq3.2}
		\int_0^T 1_{(a,b]}\d B_t^F = B_b^F-B_a^F=\int_0^T K_F^* 1_{(a,b]} \d B_s.
	\end{equation}
	For any $\psi\in \HH$, there exists a sequence $(\psi_n)_{n\in\N}\subset \mathcal S$ such that $\psi_n\to\psi$ in $\HH$, which is equivalent to $K_F^*\psi_n\to K_F^*\psi$ in $L^2([0,T];\d x)$. Thus, by Definition~\ref{defInt}, we have
	\begin{equation}
		\begin{aligned}
			\int_0^T \psi(s) \d B_s^F
			& = \lim_{n\to \infty} \int_0^T \psi_n(s) \d B_s^F
			\\ &=\lim_{n\to \infty} \int_0^T K_F^*\psi_n(s) \d B_s
			\\ &=\int_0^T K_F^*\psi(s) \d B_s
		\end{aligned}
	\end{equation}
	All limits are in $L^2(\P)$; the first uses the definition of $J^F$, the second follows from \eqref{eq3.2}, and the last follows from the It\^o isometry.
\end{proof}

\begin{remark}
	Decreusefond--\"{U}st\"{u}nel \cite{DU99} studied Wiener and divergence integrals for fractional Brownian motion. Further constructions for $H<1/2$ appear in \cite{ALN,CN}.
\end{remark}

\section{Functional inequalities (spectral gaps)}

In this section, we study the relation between spectral gaps of $\LL_{OU}$ and $\LL_{L^2}$ and the Volterra kernel $F$. We first establish quasi-invariance and construct the associated quasi-regular Dirichlet forms.

\subsection{Quasi-invariance theorem}

We prove the quasi-invariance theorem with respect to $\mu^F$.

\begin{thm}[\textbf{Quasi-invariance theorem}]
	\label{thmQuasi}
	For each $h\in\H$ and $h^F=K_F\dot{h}\in\H^F$, we have
	$$
	\int_{\Omega}G(B^F+h^F)\d\P
	=\int_{\Omega}G(B^F)\beta_{h^F}\d\P
	$$
	for any bounded Borel measurable function $G$ on $W_T$, where
\begin{equation}\label{eqAlpha}
\beta_{h^F}=
	\exp\left\{\int_0^T\dot{h}(t)\d B_t-\frac{1}{2}\|h\|_{\H}^2\right\}
\end{equation}
\end{thm}

\begin{proof}
Put $u=\dot h$ and $Z_u=\int_0^Tu(s)\d B_s$. For any times $t_1,\ldots,t_m$, the vector $(B_{t_1}^F,\ldots,B_{t_m}^F,Z_u)$ is centered Gaussian, with
\[
\operatorname{Cov}(B_{t_i}^F,Z_u)=(K_Fu)(t_i)=h^F(t_i).
\]
The Gaussian exponential formula therefore shows that exponential tilting by $\exp(Z_u-\|u\|_2^2/2)$ shifts the mean of this vector of path evaluations by $(h^F(t_i))_{i=1}^m$ and leaves its covariance unchanged. Thus the asserted identity holds for bounded evaluation cylinder functions. Since path evaluations generate the Borel sigma-field of $W_T$, a monotone class argument proves the result.
\end{proof}

By the above quasi-invariance theorem, we get the following the integration by parts formula.
\begin{thm}[\textbf{Integration by parts formula}]
\label{thmParts}
For every $h^F=K_F\dot{h}\in\H^F, h\in \H$ and $G\in\F C_b^{L^2}$,
\[
\int_\Omega D_{h^F}G(B^F)\d\P=\int_\Omega G(B^F)\Big(\int_0^T\dot{h}(s)\d B_s\Big)\d\P.
\]

\end{thm}
	
\begin{proof}
	For any $\varepsilon\in[-1,1]$ and $h^F\in \H^F$, we have $\varepsilon h^F\in \H^F$.
	By Theorem \ref{thmQuasi}, we have
	\begin{equation}\label{eq4.5}
		\int_{\Omega}G(B^F+\varepsilon h^F)\d\P=\int_{\Omega} G(B^F)\beta_{\varepsilon h^F}\d\P.
	\end{equation}
	Hence, by the dominated convergence theorem, we may differentiate both sides of \eqref{eq4.5} with respect to $\varepsilon$. Moreover,
	\begin{align*}
		\frac{\d\beta_{\varepsilon h^F}}{\d\varepsilon}\Big|_{\varepsilon=0}
		&=\frac{\d}{\d\varepsilon}\left.\exp\left\{\varepsilon\int_0^T\dot{h}(s)\d B_s-\tfrac{\varepsilon^2}{2}\|h^F\|_{\H^F}^2\right\}\right|_{\varepsilon=0}\\
		&=\int_0^T\dot{h}(s)\d B_s.
	\end{align*}
	This proves the formula.
\end{proof}

\begin{remark}
	When $B^F$ is the fractional Brownian motion, Decreusefond--\"Ust\"unel \cite{DU99} proved the quasi-invariance theorem.
\end{remark}

\subsection{Quasi-regular Dirichlet forms}

In the following, we will construct three Quasi-regular Dirichlet forms on path space. To do that, we first give several gradient definitions.

\begin{defn}\label{def4.4}
	$(1)$ \textbf{Directional derivative:} For each $h\in L^2([0,T];\d x)$ and any function $G:W_T\rightarrow\R$, we say that $D_hG$ is the directional derivative with respect to $h$ if the limit
	\begin{equation}\notag
		D_{h}G(\gamma)=\lim\limits_{\varepsilon\rightarrow0}\frac{G(\gamma+\varepsilon h)-G(\gamma)}{\varepsilon}
	\end{equation}
	exists in $L^2(\mu^T)$.
	
	$(2)$ \textbf{Gradient:} Suppose that the directional derivative $D_{h}G$ exists. By the Riesz representation theorem, there exist three gradients $\nabla^F G$, $\nabla G$ and $D G$ such that
\begin{equation}\notag
		\begin{aligned}
		&D_{h}G(\gamma)=\<\nabla^F G(\gamma), h\>_{\H^F},\quad ~~~~\forall~ h\in\H^F,\\
		&D_{h}G(\gamma)=\<\nabla G(\gamma), h\>_{\H},\quad ~~~~~~~~\forall~h\in\H,\\
		&D_{h}G(\gamma)=\<D G(\gamma), h\>_{L^2([0,T];\d x)},\quad \forall~ h\in L^2([0,T];\d x).
		\end{aligned}
	\end{equation}
	Here $\nabla^FG(\gamma)$, $\nabla G(\gamma)$ and $DG(\gamma)$ are called the damped gradient, Malliavin gradient and the $L^2$-gradient of $G$ respectively.
\end{defn}

\begin{remark}\label{rem4.6}

$(1)$
The class of Ornstein--Uhlenbeck cylinder functions is defined by
\begin{equation}\notag
	\F C_b^{OU}:=\left\{G(\gamma)=g\bigl(\gamma_{t_1},\ldots,\gamma_{t_m}\bigr):
	m\geq1,\ 0<t_1<\cdots<t_m\leq T,\ g\in C_b^1(\R^m)\right\}.
\end{equation}
As shown in \cite{CWZZ,RWZZ}, $\F C_b^{OU}$ and $\F C_b^{L^2}$ are equivalent cores for the OU Dirichlet form.
Therefore, we only need to verify the following relations for $G\in \F C_b^{OU}$ with
$$G(\gamma)=g\left(\gamma_{t_1},\cdots,\gamma_{t_m}\right),\quad\gamma\in W_T.
$$
Then for each $h^F=K_F \dot{h}\in \H^F$,
$$\begin{aligned}
D_{h^F}G(\gamma)&=\sum_{i=1}^{m}\partial_ig(\gamma)h^F(t_i)=\sum_{i=1}^{m}\partial_ig(\gamma)\<K_F^*1_{[0,t_i]},h'\>_{L^2([0,T];\d x)}\\
&=\left\<\sum_{i=1}^m\partial_igK_FK_F^*1_{[0,t_i]},h^F\right\>_{\H^F}=\left\<K_FK_F^*\left(\sum_{i=1}^m1_{[0,t_i]}\partial_ig\right),h^F\right\>_{\H^F},
\end{aligned}
$$
where $\partial_ig$ is the derivative of the $i$-th component of the function $g$.
Thus by definition \ref{def4.4},
\begin{equation}\label{eqNablaF}
	\nabla^F G(\gamma)=K_FK_F^*\left(\sum_{i=1}^m1_{[0,t_i]}\partial_ig\right)=K_FK_F^*(\dot{\overbrace{\nabla G}}(\gamma)),
\end{equation}
where
$$\dot{\overbrace{\nabla G}}(t):=\frac{d}{dt}\nabla G(t)=\sum_{i\leq m} 1_{\{t\leq t_i\}} \pp_i g.
$$

$(2)$ Suppose that $G\in \F C_b^{L^2}$ with
$$G(\gamma)=g\bigg(\int_0^Tg_1(s,\gamma_s)\d s,\dots,\int_0^Tg_m(s,\gamma_s)\d s\bigg),\quad\gamma\in W_T
$$
for $g\in C^1_b(\R^m),g_i\in C^{0,1}([0,T]\times\R;\R),1\le i\le m$. Then for each $h\in L^2([0,T];\d x)$,
$$D_hG(\gamma)=\sum_{i=1}^{m}\partial_ig\int_0^T g'_{i,2}(s,\gamma_s)h(s)\d s=\left\<\sum_{i=1}^{m}\partial_ig\, g'_{i,2}(s,\gamma_s),h(s)\right\>_{L^2([0,T];\d x)}.
$$
Thus by the definition \ref{def4.4},
\begin{equation}\label{eq4.16}
	DG(\gamma)(t)=\sum_{i=1}^{m}\partial_ig\, g'_{i,2}(t,\gamma_t).
\end{equation}
\end{remark}

\begin{defn}
We define three quadratic forms:

$(1)$ The damped quadratic form $(\E,\F C^{L^2}_b)$
		\begin{equation}\notag
			\E(\psi,\psi):=\int_{W_T}\left\|\nabla^F\psi\right\|^2_{\H^F}\d\mu^F,\quad\psi\in\F C^{L^2}_b;
		\end{equation}

$(2)$ The O-U quadratic form $(\E_{OU},\F C_b^{L^2})$
		\begin{equation}\notag
			\E_{OU}(\phi,\phi):=\int_{W_T}\left\|\nabla\phi\right\|^2_{\H}\d\mu^F,\quad\phi\in\F C_b^{L^2};
		\end{equation}

$(3)$ The $L^2$ quadratic form $(\E_{L^2},\F C_b^{L^2})$
		\begin{equation}\notag
			\E_{L^2}(\varphi,\varphi):=\int_{W_T}\left\|D\varphi\right\|^2_{L^2([0,T];\d x)}\d \mu^F, \quad\varphi\in\F C^{L^2}_b.
		\end{equation}
\end{defn}

\begin{thm}\label{thmDirichlet}
$(1)$ The damped quadratic form $(\E, \F C_b^{L^2})$ is closable, and its closure $(\E,\D(\E))$ is a quasi-regular Dirichlet form;

$(2)$  The O-U quadratic form $(\E_{OU},\F C_b^{L^2})$ is closable if and only if $\overline{\H\cap\H^F}^{\,\H}=\H$. If it is closable, then its closure is quasi-regular;

$(3)$ The $L^2$ quadratic form $(\E_{L^2},\F C_b^{L^2})$ is closable if and only if $\overline{\H^F}^{\,L^2}
		=L^2([0,T];\d x)$. If it is closable, then the corresponding closure is quasi-regular.
\end{thm}

\begin{remark}
$(1)$ When the based space is a compact Riemannian manifold, Driver--R\"ockner \cite{DR} obtained quasi-regular O-U Dirichlet forms for the Wiener measure on the associated path space; L\"obus \cite{Lob04} constructed a class of processes on the path space over a compact Riemannian manifold with unbounded diffusion; Wang--Wu \cite{WW2008,WW2009} extended Driver--R\"ockner's result to a noncompact Riemannian manifold with the integrable condition for the Ricci curvature; and Chen--Wu \cite{ChenWu14} further showed that the conclusions hold for a stochastically complete Riemannian manifold without any curvature conditions.

$(2)$ R\"ockner--Wu--Zhu--Zhu \cite{RWZZ} established quasi-regular $L^2$ Dirichlet forms for the Wiener measure on the path space over a compact manifold; Chen--Wu--Zhu--Zhu \cite{CWZZ} extended this to a noncompact Riemannian manifold.

$(3)$ For $F$-Gaussian processes, Chen--Sun--Wu \cite{CSW} obtained similar conclusions under stronger conditions.
\end{remark}

\begin{proof}[{\bf Proof of Theorem \ref{thmDirichlet}}]
The proof and notations mainly refer to \cite{CWZZ} and \cite{CSW} (see also \cite{MM,RWZZ,WW2008,WW2009}).

\textbf{Closability of the damped form}: Let $\{\psi_{n}\}_{n\geq 1} \subset \F C_{b}^{L^2}$ be a sequence of cylinder functions such that
$$\lim _{n \to \infty} \int_{\Omega} \psi_{n}^{2}\left(B^{F}\right) \d \P=0 $$
and
$$
\lim _{n,m \to \infty} \E\left(\psi_{n}-\psi_{m}, \psi_{n}-\psi_{m}\right)=\lim _{n,m \to \infty} \int_{\Omega}\left\| \nabla^{F} \psi_{n}-\nabla^{F} \psi_{m}\right\| _{\mathbb{H}^{F}}^{2}\left(B^{F}\right) \d \P=0 .$$
Then ${\nabla^{F} \psi_{n} \circ B^{F}}$ is a Cauchy sequence in $L^{2}(\Omega \to \mathbb{H}^{F} ; \mathbb{P})$ . Thus there exists a limit $\Psi \in L^{2}(\Omega \to \mathbb{H}^{F} ; \mathbb{P})$ such that
\[
\lim _{n \to \infty} \nabla^{F} \psi_{n} \circ B^{F}=\Psi.
\]
In order to prove closability, it suffices to show that $\Psi=0$. Let $\{h_{k}^{F}\}_{k\geq 1} \subset \mathbb{H}^{F}$ be an orthonormal basis in $\mathbb{H}^{F}$, and for each $k$ let $u_k\in M$ be the minimal preimage of $h_k^F$, i.e. $K_Fu_k=h_k^F$. By Theorem \ref{thmParts}, for any $G \in \F C_{b}^{L^2}$, we have
\begin{equation}\label{eq4.7}
	\begin{aligned} & \int_{\Omega} G\left(B^{F}\right) D_{h_{k}^{F}} \psi_{n}\left(B^{F}\right) \d \P \\ & =\int_{\Omega} D_{h_{k}^{F}}\left(\psi_{n}  G\right)\left(B^{F}\right) \d \P-\int_{\Omega} \psi_{n}\left(B^{F}\right) D_{h_{k}^{F}} G\left(B^{F}\right) \d \P
\\ & =\int_{\Omega} \psi_{n}\left(B^{F}\right)\left(G\left(B^{F}\right)\left(\int_0^T u_k(s)\d B_s\right)-\left<\nabla^{F} G\left(B^{F}\right), h_{k}^{F}\right>_{\mathbb{H}^{F}}\right) \d \P . \end{aligned}
\end{equation}
Since $G(B^{F})$, $\nabla^{F} G(B^{F})$ and $\int_0^T u_k(s)\d B_s$ are $L^{2}$ -integrable, by the dominated convergence theorem and \eqref{eq4.7} we get
\[\begin{aligned} & \int_{W_{T}}\left<\Psi, h_{k}^{F}\right>_{\mathbb{H}^{F}} G \d \mu^F =\lim _{n \to \infty} \int_{\Omega}\left<\nabla^{F} \psi_{n}\left(B^{F}\right), h_{k}^{F}\right>_{\mathbb{H}^{F}} G\left(B^{F}\right) \d \P \\
	& =\lim _{n \to \infty} \int_{\Omega} G\left(B^{F}\right) D_{h_{k}^{F}} \psi_{n}\left(B^{F}\right) \d \P \\
	& =\lim _{n \to \infty} \int_{\Omega} \psi_{n}\left(B^{F}\right)\left(G\left(B^{F}\right)\left(\int_0^T u_k(s)\d B_s\right)-\left<\nabla^{F} G\left(B^{F}\right), h_{k}^{F}\right>_{\mathbb{H}^{F}}\right) \d \P \\
	& =0 \end{aligned}\]
for any $G \in \F C_{b}^{L^2}$ and $k \geq 1$. This implies that $\Psi \equiv 0$.

\textbf{Closability of the O-U form}:
\emph{sufficiency.} Assume that
	$\overline{\H\cap\H^F}^{\,\H}=\H$.
Same setting with $\E_{OU}$ like above assume that there exists $\Psi\in L^2(\mu^F;\H)$ such that
\[
	\nabla G_n\longrightarrow \Psi
	\quad\text{in }L^2(\mu^F;\H).
\]
It remains to show that $\Psi=0$. For every $h\in\H^F$ and
$\Phi\in\F C_b^{L^2}$, letting $u\in M$ be the minimal preimage of $h$, Theorem \ref{thmParts} and the product rule give
\begin{equation}\label{eq:Close-OU}
\begin{aligned}
	\int_{W_T}\Phi D_hG_n\d\mu^F
	&=\int_{W_T}D_h(\Phi G_n)\d\mu^F
	  -\int_{W_T}G_nD_h\Phi\d\mu^F\\
	&=\int_{\Omega}G_n(B^F)
	  \Bigl(\Phi(B^F)\int_0^T u(s)\d B_s-D_h\Phi(B^F)\Bigr)\d\P.
\end{aligned}
\end{equation}
By the condition $\overline{\H\cap\H^F}^{\,\H}=\H$, we may choose an orthonormal basis
$\{h_k\}_{k\geq1}\subset\H\cap\H^F$ of $\H$. For each $k$, let
$u_k\in M$ be the minimal preimage of $h_k$ regarded as an element of
$\H^F$. According to \eqref{eq:Close-OU} and $\Phi\int_0^T u_k(s)\d B_s-D_{h_k}\Phi\in L^2(\mu^F)$, we have
\[
\begin{aligned}
	\int_{W_T}\Phi\langle \Psi,h_k\rangle_{\H}\d\mu^F
	&=\lim_{n\to\infty}\int_{W_T}\Phi
	  \langle\nabla G_n,h_k\rangle_{\H}\d\mu^F\\
	&=\lim_{n\to\infty}\int_{\Omega}G_n(B^F)
	  \Bigl(\Phi(B^F)\int_0^T u_k(s)\d B_s-D_{h_k}\Phi(B^F)\Bigr)\d\P=0.
\end{aligned}
\]
The arbitrary of $\Phi$ implies that
$\langle \Psi,h_k\rangle_{\H}=0$ for every $k\geq1$. Hence $\Psi=0$, and
therefore $(\E_{OU},\F C_b^{L^2})$ is closable.

\emph{Necessity.} Suppose that $\E_{OU}$ is closable. For $a\in C([0,T])$, set
\[
	\ell_a(\gamma):=\int_0^Ta(s)\gamma_s\d s,
\]
and let $u_a\in\H^F$ be the Riesz representer of $\ell_a|_{\H^F}$, i.e. $u_a=K_F\K_Fa$.
Consider the operator
\[
\begin{aligned}
V:&~\{u_a:a\in C([0,T])\}\subset\H^F&\longrightarrow\H,\\
(Vu_a)(t)&:=\int_0^t\int_r^Ta(s)\d s\d r.
\end{aligned}
\]
Fubini's theorem and the reproducing-kernel identity give
\begin{equation}\label{eq:linear-representers}
	\ell_a(h)=\langle Vu_a,h\rangle_{\H},\qquad
	\ell_a(k)=\langle u_a,k\rangle_{\H^F},\qquad
	\|u_a\|_{\H^F}^2=\operatorname{Var}_{\mu^F}(\ell_a)
\end{equation}
for $h\in\H$ and $k\in\H^F$.

The operator $V$ is well defined: if $u_a=0$, then $\ell_a=0$ $\mu^F$-a.s., so
$G=\sin(\ell_a)=0$ in $L^2(\mu^F)$, whereas $\nabla G=Vu_a$ $\mu^F$-a.s.; closability therefore gives $Vu_a=0$.
Its domain is dense in $\H^F$, since \eqref{eq:linear-representers} and the density of $C([0,T])$ in $L^2([0,T];\d x)$ show that a vector orthogonal to every $u_a$ must vanish.

Moreover, $V$ is closable. Indeed, suppose that $u_{a_n}\to0$ in $\H^F$ and $Vu_{a_n}\to h$ in $\H$. By \eqref{eq:linear-representers}, $\ell_{a_n}\to0$ in $L^2(\mu^F)$. Hence
\[
	G_n:=\sin(\ell_{a_n})\longrightarrow0,
	\qquad
	\nabla G_n=\cos(\ell_{a_n})Vu_{a_n}\longrightarrow h
\]
in the corresponding $L^2$ spaces; here the second convergence follows from
$\cos(\ell_{a_n})\to1$ in $L^2(\mu^F)$. The closability of $\E_{OU}$ implies $h=0$.

It remains to compute the adjoint domain. By \eqref{eq:linear-representers},
\[
\begin{aligned}
	h\in\Dom(V^*)
	&\Longleftrightarrow
	\exists\,k\in\H^F:\
	\ell_a(h)=\ell_a(k)\quad\forall a\in C([0,T])\\
	&\Longleftrightarrow h=k\in\H\cap\H^F.
\end{aligned}
\]
Thus $\Dom(V^*)=\H\cap\H^F$. Since a densely defined closable operator has a densely defined adjoint, we conclude that
\[
	\overline{\H\cap\H^F}^{\,\H}=\H,
\]
as required.

\textbf{Closability of the $L^2$ form}:
The sufficiency proof is the same as for the O-U form.

\emph{Necessity.} Suppose instead that $\H^F$ is not dense in $L^2([0,T];\d x)$. Choose a non-zero $q\in L^2([0,T];\d x)$ such that $q\perp\H^F$. Since $\H^F=\operatorname{Ran}K_F$, Hilbert-space duality gives
\begin{equation}\label{eq:adjoint-kernel-q}
	\K_F^*q=0.
\end{equation}
For $X_q(\gamma):=\int_0^Tq(t)\gamma_t\d t$, stochastic Fubini's theorem gives
\[
	X_q(B^F)=\int_0^T(\K_F^*q)(s)\d B_s=0
	\qquad\P\text{-a.s.}
\]
Choose $q_n\in C([0,T])$ with $q_n\to q$ in $L^2([0,T];\d x)$, and set
\[
	X_n(\gamma)=\int_0^Tq_n(t)\gamma_t\d t,
	\qquad G_n(\gamma)=\sin(X_n(\gamma)).
\]
Since $\K_F^*$ is bounded on $L^2([0,T];\d x)$,
\[
	\|X_n-X_q\|_{L^2(\mu^F)}
	=\|\K_F^*(q_n-q)\|_{L^2([0,T];\d x)}\longrightarrow0.
\]
Thus $G_n\to0$ in $L^2(\mu^F)$. Passing to a subsequence, assume also that $X_n\to0$ $\mu^F$-a.s. The explicit $L^2$-gradient formula gives
\[
	DG_n=\cos(X_n)q_n.
\]
Consequently,
\[
\begin{aligned}
	\EE_{\mu^F}\|DG_n-q\|_{L^2}^2
	&\leq2\|q_n-q\|_{L^2}^2
	+2\|q\|_{L^2}^2\EE_{\mu^F}|\cos(X_n)-1|^2\longrightarrow0,
\end{aligned}
\]
by dominated convergence. Hence $DG_n\to q\ne0$ in
$L^2(\mu^F;L^2([0,T];\d x))$, contradicting closability.

\textbf{Quasi-regularity}

$(a)$ \emph{The case of damped Dirichlet form}: It is obvious that $\F C_{b}^{L^2}$ is dense in $\D(\E)$ under the $\E_1^{1/2}$-norm, and has a countable dense subset separating points on $W_T$. Then, by \cite[Definition IV-3.1]{MM}, to verify the quasi-regularity of $(\E,\D(\E))$ it suffices to find a sequence of compact sets $\{K_n\}\subset W_T$ such that
\[\lim _{n \to \infty} \Cap\left(W_{T} \backslash K_{n}\right)=0,\]
where Cap is the capacity induced by $(\E, \D(\E))$.

$(I)$
We intend to prove that for some constant $C>0$ and any $\gamma \in W_{T}$ , the function
\[F_{\gamma }(\sigma ):=\sup_{t\in [0,T]}|\gamma _{t}-\sigma _{t}| ,\quad  \sigma \in W_{T}\]
is in $\D(\E)$ and
\begin{equation}
	\label{eq4.8}
	\E (F_{\gamma },F_{\gamma })\leq C.
\end{equation}
Let $\{t_{k}\}_{k=1}^{\infty}$ be a countable dense subset of $[0, T]$.
For each $\gamma \in W_{T}$ , define
\[F_{\gamma ,k}(\sigma ):=|\gamma (t_{k})-\sigma (t_{k})|,
\quad F_{\gamma }^{(n)}:=\operatorname* {max}_{1\leq k\leq n}F_{\gamma ,k},
\quad k,n\in \N,\quad \sigma \in W_{T}.\]
According to \eqref{eqNablaF} we derive
\[\left\| \nabla^{F} F_{\gamma, k}\right\| _{\H^{F}}^{2}=\left\| K_F^* \dot{\overbrace{\nabla F_{\gamma, k}}}\right\| _{L^2([0,T];\d x)}^{2}=\|K_F^*1_{[0,t_k]}\|_{L^2([0,T];\d x)}^2 = \int_0^T F(t_k,s)^2\d s.\]
Then for each $k\geq1$,
\[\E(F_{\gamma,k},F_{\gamma,k})=\int_0^T F(t_k,s)^2\d s.\]
Thus, for any $n\geq1$,
\[\E(F_{\gamma}^{(n)},F_{\gamma}^{(n)})\leq\max_{1\leq k\leq n}\E(F_{\gamma,k},F_{\gamma,k})\leq\max_{1\leq k\leq n}\int_0^T F(t_k,s)^2\d s.\]
Since $F_{\gamma}^{(n)}\uparrow F_\gamma$ as $n\uparrow\infty$ and $F_\gamma\in L^2(W_T;\mu^F)$, it follows from \cite[Lemma I-2.12]{MM} that $F_\gamma\in\D(\E)$ and \eqref{eq4.8} holds with $C=\sup_{t\in[0,T]}\int_0^T F(t,s)^2\d s$.

$(II)$ Now let ${\gamma_{k}: k \geq 1}$ be a countable dense subset of $W_{T}$ . For any $n \geq 1$ , let
\[F_{n}:=\inf _{1 \leq k \leq n} F_{\gamma_{k}} .\]
Similar to the above argument, we obtain $F_n\in\D(\E)$ and
\[\E(F_n,F_n)\leq C,\quad n\geq1.\]
Since $0\leq F_{n+1}\leq F_n$ and $F_1\in L^2(W_T;\mu^F)$, we obtain
\[\|F_n\|_{L^2(W_T;\mu^F)}\leq\|F_1\|_{L^2(W_T;\mu^F)},\quad n\geq1.\]
Therefore,
\[\E_1(F_n,F_n)\leq C+\|F_1\|_{L^2(W_T;\mu^F)}^2=:C_2,\quad n\geq1.\]
Since $\{\gamma_i:i\geq1\}$ is dense in $W_T$, we have $F_n\downarrow0$ pointwise as $n\uparrow\infty$.

$(III)$ Now, following the line of \cite{DR} (see pages 606 and 607 therein), we need to find a sequence of closed subsets $\{K_m\subset W_T:m\in\mathbb N\}$ such that $F_n\to0$ uniformly on every $K_m$ and the corresponding capacity of $W_T\setminus K_m$ tends to zero. Since $\{F_n\}$ is uniformly bounded in the form norm, the Banach--Saks argument gives a subsequence, denoted again by $\{F_n\}$, such that
\[\E(\bar F_n,\bar F_n)\to0,\quad n\uparrow\infty,\]
where $\bar F_n:=n^{-1}\sum_{i=1}^nF_i$. According to \cite[Lemma I-2.12 and Propositions \uppercase\expandafter{\romannumeral3}-3.5]{MM}, there exist a further subsequence and an $\E$-nest $\{K_m\}$ such that $\bar F_n$ converges uniformly to zero on every $K_m$. Since $F_n$ is decreasing, $F_n\leq\bar F_n$, and hence $F_n\to0$ uniformly on every $K_m$ as well. As in \cite{RS}, proof of Proposition 3.1, or in \cite{DR}, proof of Proposition 5, it follows now from the definition of $F_{n}$ that every $K_m$ is totally bounded. Thus, $\{K_m\}$ is an $\E$-nest consisting of compact sets, which proves the quasi-regularity of the damped form on $W_T$.

$(b)$ \emph{The case of O-U Dirichlet form}: The same procedure yields quasi-regularity of the O-U form on $W_T$. Indeed, the Cameron--Martin representer of evaluation at $t_k$ is $s\mapsto s\wedge t_k$, and hence
\[
	\|\nabla F_{\gamma,k}\|_{\H}^2=t_k\leq T.
\]
Therefore, the preceding argument applies with $\E_{OU}$ in place of $\E$ and with $T$ replacing $C$.

$(c)$ \emph{The case of $L^2$ Dirichlet form}: Put
\[
	J:W_T\longrightarrow L^2([0,T];\d x),\quad J(\gamma)=[\gamma],
	\qquad \nu:=J_\#\mu^F.
\]
Here $[\gamma]$ denotes the $L^2([0,T];\d x)$-equivalence class represented by the continuous path $\gamma$, and $J_\#\mu^F$ denotes the pushforward of $\mu^F$ under $J$.
We regard the closure, still denoted by $\E_{L^2}$, as a form on $L^2(L^2([0,T];\d x);\nu)$, $\nu$ is the extension of $\mu$ in $L^2([0,T];\d x)$.
According to the closability condition for $\E_{L^2}$ and $(2)$, we know that
\[
	\supp\nu=\overline{\operatorname{Ran}K_F}^{\,L^2([0,T];\d x)}
	=\overline{\H^F}^{\,L^2([0,T];\d x)}=L^2([0,T];\d x).
\]

We first establish the analogue of \eqref{eq4.8}. For $\gamma\in L^2([0,T];\d x)$, let
\[
	F_\gamma(x):=\|x-\gamma\|_{L^2([0,T];\d x)},\qquad x\in L^2([0,T];\d x).
\]
Choose an orthonormal basis $\{e_j\}_{j\geq1}\subset C([0,T])$ of $L^2([0,T];\d x)$, and let $P_n$ be the orthogonal projection onto $\operatorname{span}\{e_1,\ldots,e_n\}$. For $\varepsilon>0$, set
\[
	F_\gamma^{(n,\varepsilon)}(x)
	:=\sqrt{\|P_n(x-\gamma)\|_{L^2([0,T];\d x)}^2+\varepsilon^2}-\varepsilon.
\]
Its $L^2$-gradient is
\[
	\frac{P_n(x-\gamma)}
	{\sqrt{\|P_n(x-\gamma)\|_{L^2([0,T];\d x)}^2+\varepsilon^2}},
\]
whose $L^2([0,T];\d x)$-norm is bounded by $1$.
To justify that these functions belong to the closed domain, choose $\chi_R\in C_b^\infty(\R)$ such that $\chi_R(r)=r$ for $0\leq r\leq R$, $\sup_R\|\chi_R'\|_\infty<\infty$, and $\chi_R(r)\to r$, $\chi_R'(r)\to1$ for every $r\geq0$. Then
\[
	G_\gamma^{(n,\varepsilon,R)}
	:=\chi_R\bigl(F_\gamma^{(n,\varepsilon)}\bigr)
	\in\F C_b^{L^2}.
\]
Indeed, its coordinates are
$\langle x,e_j\rangle_{L^2([0,T];\d x)}=\int_0^Te_j(s)x(s)\d s$, which are allowed in the definition of $\F C_b^{L^2}$. The chain rule and dominated convergence yield, as $R\to\infty$,
\[
	G_\gamma^{(n,\varepsilon,R)}\longrightarrow
	F_\gamma^{(n,\varepsilon)}
	\quad\text{in }L^2(L^2([0,T];\d x);\nu),
\]
while their gradients converge in $L^2(L^2([0,T]);L^2([0,T]); \nu)$. Thus they converge in the $\E_{L^2,1}^{1/2}$-norm, and
$F_\gamma^{(n,\varepsilon)}\in\D(\E_{L^2})$.

Take $\varepsilon_n\downarrow0$. Since $\nu$ has finite second moment, dominated convergence gives
\[
	F_\gamma^{(n,\varepsilon_n)}\longrightarrow F_\gamma
	\quad\text{in }L^2(L^2([0,T];\d x);\nu),
\]
while the corresponding gradient fields converge to
\[
	\begin{cases}
	\dfrac{x-\gamma}{\|x-\gamma\|_{L^2([0,T];\d x)}},&x\neq\gamma,\\[1.2ex]
	0,&x=\gamma,
	\end{cases}
	\quad\text{in }L^2(L^2([0,T]);L^2([0,T]); \nu).
\]
Consequently, $\{F_\gamma^{(n,\varepsilon_n)}\}_{n\geq1}$ is Cauchy in the $\E_{L^2,1}^{1/2}$-norm. The closedness of $\E_{L^2}$ therefore implies
\begin{equation}\label{eq:L2-distance-energy}
	F_\gamma\in\D(\E_{L^2}),\qquad
	\E_{L^2}(F_\gamma,F_\gamma)\leq1,
	\qquad \gamma\in L^2([0,T];\d x).
\end{equation}

Let $\{\gamma_k:k\geq1\}$ be a countable dense subset of $L^2([0,T];\d x)$ and define
\[
	F_n(x):=\min_{1\leq k\leq n}\|x-\gamma_k\|_{L^2([0,T];\d x)}.
\]
The lattice and locality properties of the closed form, together with \eqref{eq:L2-distance-energy}, show that
\[
	F_n\in\D(\E_{L^2}),\qquad
	\E_{L^2}(F_n,F_n)\leq1.
\]
Furthermore,
\[
	0\leq F_{n+1}\leq F_n\leq F_1,
	\qquad F_n(x)\downarrow0\quad(x\in L^2([0,T];\d x)).
\]
Consequently $F_n\to0$ in $L^2(L^2([0,T];\d x);\nu)$ and
\[
	\sup_{n\geq1}\E_{L^2,1}(F_n,F_n)
	\leq1+\|F_1\|_{L^2(L^2([0,T];\d x);\nu)}^2<\infty.
\]
Thus $F_n$ converges weakly to zero in the Hilbert space
$(\D(\E_{L^2}),\E_{L^2,1})$. By the Banach--Saks theorem, there is a subsequence $\{F_{n_j}\}$ such that
\[
	\overline F_N:=\frac1N\sum_{j=1}^NF_{n_j},
	\qquad
	\E_{L^2,1}(\overline F_N,\overline F_N)\longrightarrow0.
\]
By \cite[Lemma I-2.12 and Proposition \uppercase\expandafter{\romannumeral3}-3.5]{MM}, a further subsequence admits an $\E_{L^2}$-nest of closed sets $\{K_m\}$ on each of which $\overline F_N$ converges uniformly to zero. Since $F_n$ is decreasing,
\[
	F_{n_j}\leq\frac1j\sum_{i=1}^jF_{n_i}=\overline F_j,
\]
so the entire sequence $F_n$ converges uniformly to zero on every $K_m$.

Each $K_m$ is compact in $L^2([0,T];\d x)$. Indeed, for every $\varepsilon>0$, uniform convergence gives an $N$ such that
\[
	K_m\subset\bigcup_{k=1}^N B_{L^2([0,T];\d x)}(\gamma_k,\varepsilon).
\]
Thus $K_m$ is totally bounded. It is also closed, and $L^2([0,T];\d x)$ is complete, so $K_m$ is compact. Hence
\[
	\Cap_{L^2}(L^2([0,T];\d x)\setminus K_m)\longrightarrow0.
\]

Since $\D(\E_{L^2})$ is separable under the $\E_{L^2,1}^{1/2}$-norm, a countable subset of these continuous cylinder functions is dense in the form domain and hence provides the required quasi-continuous functions. Finally, if $\{h_j\}_{j\geq1}\subset C([0,T])$ is dense in $L^2([0,T];\d x)$, the continuous cylinder functions
\[
	u_j(x):=\tanh\langle x,h_j\rangle_{L^2([0,T];\d x)},\qquad j\geq1,
\]
separate the points of $L^2([0,T];\d x)$. We have therefore verified all the conditions in \cite[Definition IV-3.1]{MM}; hence $\E_{L^2}$ is quasi-regular on $(L^2([0,T];\d x),\nu)$.

\end{proof}

\begin{remark}
	\label{rmk:closability-comparison}
	The OU closable criterion implies the $L^2$ criterion. If $\overline{\H\cap\H^F}^{\,\H}=\H$, then the continuous embedding gives
	\[
		\H\subset\overline{\H^F}^{\,L^2([0,T];\d x)},
	\]
	and, together with $\overline{\H}^{\,L^2([0,T];\d x)}=L^2([0,T];\d x)$, this yields
	\[
		\overline{\H^F}^{\,L^2([0,T];\d x)}=L^2([0,T];\d x).
	\]

\end{remark}
\begin{remark}\label{closable cond}
	$(1)$
	For the fractional Brownian motion $B^{K_H}$, by \cite[Proposition 6]{DU99}, the characterization $\H^H=K_H(L^2([0,T]))=I_{0+}^{H+1/2}(L^2([0,T]))$. Since $C^\infty \cap \H \subset \H\cap\H^F$ is dense in $\H$, hence $\overline{\H\cap\H^F}^{\,\H}=\H$.

	$(2)$
	Fix constants $0\leq\mu<1$ and $1/2<\nu<2$. we assume that $f^{-1}\in L^2([0,T];\d x)$ and
$$x^{\nu-\mu-1}f\in L^2([0,T];\d x) ~~\text{and}~~x^{-\mu}f\in L^{\frac2{2\nu-1}+\varepsilon}([0,T];\d x)$$
 for some $\varepsilon>0$. Setting
	\[
		F(t,s)=(t-s)^{\nu-1}t^{-\mu}f(s)+\mu f(s)\int_s^t(u-s)^{\nu-1}u^{-\mu-1}\d u.
	\]
Since
	$$\left\{K_F(g/f); g\in C_c^\infty([0,T])\right\} \subset \H^F
	$$
	is dense in $\H$,
	then $\overline{\H\cap\H^F}^{\,\H}=\H$.

	$(3)$
	Assume that $\sigma\in L^{2+\varepsilon}([0,T];\d x)$ for some $\varepsilon>0$ and $\sigma\neq0~a.e.$. Let $f\in L^2([0,T];dx)$ satisfying with $f^{-1}, f'\in L^2([0,T];dx)$.
	For diffusion process $\d X_t=\frac{f'(t)}{f(t)}X_t\d t+f(t)\sigma(t) \d B_t$ starting from $X_0=0$, it is F-Gaussian process and can be expressed as
	$$
	X_t=\int_0^tF(t,s)\d B_s = \int_0^t f(t)\sigma(s)\d B_s.
	$$
	Define
\[
	A:=\left\{g(t)=f(t)\int_0^t\sigma(s)u(s)\d s\big|u\in L^2,\ \sigma u\in L^2\right\}.
	\]
Then $A\subset\H\cap\H^F$ and $A$
	is dense in $\H$. Indeed, for any $h\in \H$, taking $q=h/f, u_n=\frac{q'}{\sigma}1_{\{\sigma\geq 1/n\}}$. Then
  $g_n(\cdot):=f(\cdot)\int_0^\cdot\sigma(s)u_n(s)\d s\in A$ and $\sigma u_n\rightarrow q'$, thus
	\[\aligned
	\left\|g_n-h\right\|^2_\H&=\int^T_0\Big|f(t)\sigma(t)u_n(t)+f'(t)\int_0^t\sigma(s)u_n(s)\d s-h'(t)\Big|\d s\\
&=\int^T_0\Big|f(t)\sigma(t)u_n(t)+f'(t)\int_0^t\sigma(s)u_n(s)\d s-f(t)q'(t)-f'(t)q(t)\Big|\d s\\&
	\leq(\sqrt T\|f'\|_{L^2}+\|f\|_\infty)\|\sigma u_n-\dot q\|_{L^2}\longrightarrow0,~~\text{as}~~n\rightarrow\infty.
	\endaligned\]
	Hence $\overline{\H\cap\H^F}^{\,\H}=\H$.
\end{remark}

\subsection{Logarithmic Sobolev inequalities}

In this section, we will establish Logarithmic Sobolev inequalities with respect to Dirichlet forms $(\E,\D(\E))$, $(\E_{OU},\D(\E_{OU}))$ and $(\E_{L^2},\D(\E_{L^2}))$ respectively. To do that, we assume further in this section that these Dirichlet forms are closable.

In the following, by using the Clark-Ocone-Haussmann formula for Brownian motion $B_t$, we will derive the following damped Logarithmic Sobolev inequality for the Gaussian Wiener measure $\mu^F$.

\begin{thm}[\textbf{Damped Logarithmic Sobolev inequality}]\label{thmDamped}
	\begin{equation*}
		\Ent_{\mu^F}(G^2):=\EE_{\mu^F}\left(G^2\log\frac{G^2}{\EE_{\mu^F}(G^2)}\right)\le 2\E(G,G),\quad G\in\D(\E).
	\end{equation*}
\end{thm}

\begin{proof}
	We first work with bounded cylinder functions and apply the argument to $G^2+\varepsilon$, $\varepsilon>0$, before letting $\varepsilon\downarrow0$. For the martingale calculation below write $G$ for a bounded cylinder function bounded away from zero.
		By the classical Clark-Ocone formula of the Brownian motion $B_t$, we know that
	\begin{equation}\label{eq5.2}
		G\left(B\right)=\EE\big[G\left(B\right)\big]+\int_0^T\EE\big[\dot{\overbrace{\nabla G}}(B)(t)|\F_t\big]\d B_t,
	\end{equation}
	According to the definition of $B^F$, we know that $B^F$ may be regarded as a functional of $B$, i.e. $B^F_\cdot=\int^\cdot_0F(\cdot,s)\d B_s:=\Phi\circ B$. Now replace $G$ by $G\circ\Phi$ in \eqref{eq5.2}, we have
	\begin{equation}\notag
		G\circ\Phi\left(B\right)=\EE\big[(G\circ\Phi)\left(B\right)\big]+\int_0^T\EE\big[ \dot{\overbrace{\nabla (G\circ\Phi)}} (B)(t) |\F_t\big]\d B_t.
	\end{equation}
		Notice that for each $h\in \H$,
	$$\begin{aligned}
	\left\< \dot{\overbrace{\nabla (G\circ\Phi)}} (B),\dot{h}\right\>_{L^2([0,T];\d x)}&=D_{h}(G\circ\Phi)(B)=D_{K_F\dot{h}}G(B^F)
	\\
	&=\int_0^T  \dot{\overbrace{\nabla G}}(B^F)(t) \d K_F\dot{h}(t)\\
	&=\left\<K_F^*\dot{\overbrace{\nabla G}}(B^F),\dot{h}\right\>_{L^2([0,T];\d x)}.
	\end{aligned}
	$$
	Due to the arbitrariness of $h$, we obtain
	\begin{equation}\label{eq5.3}
		G(B^F)=
		\EE\big[G(B^F)\big]+\int_0^T\EE\big[K_F^*\dot{\overbrace{\nabla G}}(B^F)(t)|\F_t\big]\d B_t.
	\end{equation}
	Consider the right continuous version of martingale
	\begin{equation*}
		M_t:=\EE\big[G(B^F)|\F_t\big], ~~\quad 0\le t\le T,
	\end{equation*}
	then by \eqref{eq5.3} we have
	\begin{equation*}
		M_t=\EE\big[G(B^F)\big]+{\int_0^t\EE\big[(K_F^*\dot{\overbrace{\nabla G}}(B^F))(s)\mid\F_s\big]\d B_s}
	\end{equation*}
	and
	\begin{equation*}
		\d M_t=\EE\big[K_F^*\dot{\overbrace{\nabla G}}(B^F)(t)|\F_t\big]\d B_t.
	\end{equation*}
	Consider the function $f(x):=x\log x$. Then
	$$f'(x)=\log x+1,\quad f''(x)=x^{-1},\quad x>0.
	$$
	By the It\^{o}'s formula, we have
	\begin{equation*}
		\begin{aligned}
		f(M_T)-f(M_0)&=\int_0^T(\log(M_t)+1)\d M_t\\
		&+\frac{1}{2}\int_0^T\frac{\left(\EE\big[K_F^*\dot{\overbrace{\nabla G}}(B^F)(t)|\F_t\big]\right)^2}{M_t}\d t.
		\end{aligned}
	\end{equation*}
	Taking expectations on both sides in the above equation with respect to $\EE$, we get
	\begin{equation*}
		\begin{aligned}
		&\EE\big[G(B^F)\log G(B^F)\big]-\EE\big[G(B^F)\big]\log \EE\big[G(B^F)\big]\\
		&=\frac{1}{2}\EE\bigg[\int_0^T\frac{\left(\EE\big[K_F^*\dot{\overbrace{\nabla G}}(B^F)(t)|\F_t\big]\right)^2}{\EE\big[G(B^F)|\F_t\big]}\d t\bigg].
		\end{aligned}
	\end{equation*}
	By using the Cauchy-Schwarz inequality and replacing $G$ by $G^2$, we get
	$$\begin{aligned}
\left(\EE\big[K_F^*\dot{\overbrace{\nabla G^2}}(B^F)(t)|\F_t\big]\right)^2
	&=\left(\EE\big[2G(B^F) K_F^*\dot{\overbrace{\nabla G}}(B^F)(t)|\F_t\big]\right)^2\\&
	\le4\EE\big[G(B^F)^2|\F_t\big]\EE\big[(K_F^*\dot{\overbrace{\nabla G}}(B^F)(t))^2|\F_t\big].
	\end{aligned}
	$$
	Therefore
	$$\begin{aligned}
	\Ent_{\mu^F}(G^2)
	&\le2\EE\int_0^T\EE\big[(K_F^*\dot{\overbrace{\nabla G}}(B^F)(t))^2|\F_t\big]\d t\\
	&\le2\EE\int_0^T(K_F^*\dot{\overbrace{\nabla G}}(B^F)(t))^2\d t\\
	&=2\EE_{\mu^F}\left\|(\nabla^FG)\right\|_{\H^F}^2.
	\end{aligned}
	$$
\end{proof}

However, in general, the logarithmic Sobolev inequality with respect to the Ornstein-Uhlenbeck Dirichlet form $(\E_{OU}, \D(\E_{OU}))$ does not hold. In the following, we provide an example to show that the weak Poincar\'e inequality does not hold,  which implies the Poincar\'e and logarithmic Sobolev inequalities fail.

\begin{counterexample}[\textbf{Weak Poincar\'e inequality (O-U case)}]\label{counter:OU-Poincare}
If $K_F^*$ is unbounded, there exists a sequence $(h_n)\subset \mathcal{S}$ and
$$G_n(B^F)=\int_0^T h_n(t)\d B^F_t
\quad \text{and} \quad
H_n=\sin\left(\frac{G_n}{\|K_F^*h_n\|_{L^2}}\right),
$$
such that for any function
$\alpha:(0,1)\to(0,\infty)$, we have
$$
\lim_{n\rightarrow\infty}\frac{\Var_{\mu^F}(H_n)}{\alpha(r)\E_{OU}(H_n,H_n)+r\operatorname{Osc}H_n^2}
\geq \frac{1-\e^{-2}}{8r}, \quad \forall~1>r>0,
$$
which implies weak Poincar\'e inequality does not hold.
\end{counterexample}

\begin{proof}
By the definition of the Malliavin gradient, we get
$$
\nabla G_n(\gamma)(t)=\int_0^t h_n(s)\d s,
$$
and hence
$$
\E_{OU}(G_n,G_n)=\EE_{\mu^F}\|\nabla G_n\|_{\H}^2
=\left\|h_n\right\|_{L^2([0,T];\d x)}^2.
$$
In addition, since $K_F^*$ is unbounded, there exists a sequence $\{h_n\}_{n\geq 1}$ such that
$$
q_n=\frac{\|K_F^*h_n\|_{L^2([0,T];\d x)}^2}
{\|h_n\|_{L^2([0,T];\d x)}^2}\longrightarrow\infty.
$$
Since $G_n(B^F)$ is a centered Gaussian random variable and $\Var_{\mu^F}(G_n)=\|K_F^* h_n\|_{L^2([0,T];\d x)}^2$, then
 $G_n(B^F)/\|K_F^* h_n\|_{L^2([0,T];\d x)}$ has the standard normal distribution, the chain rule for Malliavin derivative gives
$$
\E_{OU}(H_n,H_n)
=\frac{\|\nabla G_n\|^2_\H }{\|K_F^* h_n\|_{L^2([0,T];\d x)}^2} \EE_{\mu^F} \cos\left(  \frac {G_n(B^F)}{\|K_F^* h_n\|_{L^2([0,T];\d x)}}   \right)
=\frac{1+\e^{-2}}{2q_n},
$$
whereas
$$
\Var_{\mu^F}(H_n)=\frac{1-\e^{-2}}2=:c_0,
\qquad
\operatorname{Osc}(H_n)^2\leq4.
$$
If a weak Poincar\'e inequality holds with the rate function $\alpha$ for $\mu^F$, i.e.
$$\Var_{\mu^F}(G)\leq \alpha(r)\E_{OU}(G,G)
		+r\,\operatorname{Osc}(G)^2,
		\quad r\in(0,1),\quad G\in\D(\E_{OU}).$$
Then we have
$$
c_0\leq \alpha(r)\frac{1+\e^{-2}}{2q_n}+4r, \quad r\in(0,1).
$$
Letting $n\to\infty$ yields $c_0\leq4r$,  this implies $c_0\leq \inf_{r\in (0,1)} 4r=0$. This is a contradiction to $c_0>0$.
\end{proof}

In particular, for the fractional Brownian motion $F=K_H$ with $H<1/2$, by Remark \ref{closable cond}, we know that the corresponding O-U Dirichlet form is quasi-regular. By the above example, $\|K_{K_H}^*\|=\infty$ implies that weak Poincar\'e inequality for $\mu^{K_H}$ does not hold.

However, if $K_F^*|_{\mathcal{S}}: L^2([0,T];\d x) \to L^2([0,T];\d x)$ is bounded, then the logarithmic Sobolev inequality with respect to $(\E_{OU},\D(\E_{OU}))$ is guaranteed. Furthermore, the logarithmic Sobolev inequality, Poincar\'e inequality and weak Poincar\'e inequality are equivalent to each other.

\begin{thm}[\textbf{Logarithmic Sobolev and Poincar\'e inequalities(O-U case)}]
	\label{thm:LSI-Poincare-OU}
	The following $(i)$--$(iv)$ are equivalent.
\end{thm}
\noindent $(i)$ $K_F^*|_{\mathcal{S}}: L^2([0,T];\d x) \to L^2([0,T];\d x)$ is bounded.

\noindent $(ii)$ There exists a constant $C_{LS}$ such that the logarithmic Sobolev inequality holds:
	\begin{equation}\label{eqLSI-OU}
		\Ent_{\mu^F}(G^2)\leq C_{LS}\E_{OU}(G,G),\quad G\in\D(\E_{OU}).
	\end{equation}

\noindent $(iii)$ There exists a constant $C_P$ such that the Poincar\'e inequality holds:
	\begin{equation}\label{eqPoincare-OU}
		\Var_{\mu^F}(G)\leq C_P\E_{OU}(G,G),\quad G\in\D(\E_{OU}).
	\end{equation}

\noindent $(iv)$ A weak Poincar\'e inequality holds: there is a function
$\alpha:(0,1)\to(0,\infty)$ such that
	\begin{equation}\label{eqWeakPoincare-OU}
		\Var_{\mu^F}(G)\leq \alpha(r)\E_{OU}(G,G)
		+r\,\operatorname{Osc}(G)^2,
		\quad r\in(0,1),\quad G\in\D(\E_{OU}),
	\end{equation}
where $\operatorname{Osc}(G):=\mathop{\rm ess\,sup}G-\mathop{\rm ess\,inf}G$.
\begin{proof} $(i)\Rightarrow (ii)$.

	For $G\in \F C^{OU}_b$, Theorem \ref{thmDamped} shows
	$$\Ent_{\mu^F}(G^2)\leq 2\E(G,G)=2\EE_{\mu^F}\|\nabla^F G\|_{\H^F}^2.
	$$
	By \eqref{eqNablaF}, we have
	$$\aligned \EE_{\mu^F}\|\nabla^F G\|_{\H^F}^2&=\EE_{\mu^F}\|K_F \circ K_F^*\dot{\overbrace{\nabla G}}\|_{\H^F}^2
=\EE_{\mu^F}\|K_F^* \dot{\overbrace{\nabla G}}\|_{L^2([0,T];\d x)}^2\\
&\leq \|K_F^*\|^2 \E_{OU}(G,G).
	\endaligned $$
	Therefore, the logarithmic Sobolev inequality \eqref{eqLSI-OU} holds with $C_{LS}=2\|K_F^*\|^2$.

 $(ii)\Rightarrow (iii)\Rightarrow (iv)$. It is obvious.

$(iv)$ implies $(i)$: if $K_F^*$ were unbounded, by the above Example \ref{counter:OU-Poincare}, no weak Poincar\'e inequality could hold. Consequently, all four statements are equivalent to each other.
\end{proof}

\begin{remark}$(1)$
	Assume that $\frac{\pp}{\pp t}F(t,s)\in L^1([0,T];\d x)$ with respect to $t$ and
	$$F(t,t) \text{ and }\int_0^t |\frac{\pp}{\pp t}F(t,s)|\d s\int_s^T |\frac{\pp}{\pp r}F(r,s)|\d r
	$$ are bounded. Then $K_F^*$ is bounded.

\begin{proof}
	Since $\frac{\pp}{\pp t}F(t,s)\in L^1([0,T];\d x)$ with respect to $t$, by definition \eqref{eqKF*} we have
	\begin{align*}
		K_F^*g(s)  &=  g(s) F(T,s) +\int_s^T (g(t) - g(s))\frac{\pp}{\pp t}F(t,s)\d t
		\\ &= g(s) F(s,s) +\int_s^T g(t)  \frac{\pp}{\pp t}F(t,s)\d t.
	\end{align*}
	Here $\|g(s)F(s,s)\|_{L^2([0,T];\d x)}\leq \|g\|_{L^2([0,T];\d x)}\sup_{s\in [0,T]} |F(s,s)|$ and
	\begin{align*}
		\int_0^T \d s\left(\int_s^T g(t)  \frac{\pp}{\pp t}F(t,s)\d t\right)&^2
		\leq \int_0^T \d s\int_s^T g^2(t)  |\frac{\pp}{\pp t}F(t,s)|\d t\int_s^T |\frac{\pp}{\pp t}F(t,s)|\d t
		\\ & =\int_0^T g^2(t)\d t\int_0^t |\frac{\pp}{\pp t}F(t,s)|\d s\int_s^T |\frac{\pp}{\pp r}F(r,s)|\d r
		\\ &=\|g\|_{L^2([0,T];\d x)}^2 \sup_{t\in [0,T]}\left(\int_0^t |\frac{\pp}{\pp t}F(t,s)|\d s\int_s^T |\frac{\pp}{\pp r}F(r,s)|\d r\right).
	\end{align*}
	Therefore
	$$\|K_F^*g\|_{L^2([0,T];\d x)}\leq \frac C2\|g\|_{L^2([0,T];\d x)}
	$$
	where $C$ is a constant.
\end{proof}

	$(2)$ For fractional Brownian motion with respect to $H>1/2$, it is easy to prove $K_F^*$ is bounded.
\end{remark}

From the damped logarithmic Sobolev inequality in Theorem \ref{thmDamped}, we also derive the following Logarithmic Sobolev inequality with respect to $L^2$-Dirichlet form $(\E_{L^2}, \D(\E_{L^2}))$.

\begin{thm}[\textbf{Logarithmic Sobolev inequality($L^2$-case)}]
	\begin{equation}\notag
		\Ent_{\mu^F}(G^2)\le C\E_{L^2}(G,G),
	\end{equation}
	for $G\in\D(\E_{L^2})$, where
	\begin{equation}\notag
		C=2\|F\|_{L^2([0,T]^2;\d x\times \d x)}^2=2\int_{[0,T]^2}F(t,s)^2\d s\d t.
	\end{equation}
\end{thm}

\begin{proof}
	Since $\F C_b^{L^2}\subset\D(\E)$, Theorem \ref{thmDamped} holds for each $G\in \F C_b^{L^2}$. For any $h\in\H^F$, there exists a  $\dot{h}\in L^2([0,T];\d x)$ such that $h^F=K_F\dot{h}$, thus we have
	\begin{equation}\notag
		\begin{aligned}
			\left\|\nabla^FG\right\|_{\H^F}^2&=\sup_{\|h^F\|_{\H^F}\le1}\left\<\nabla^FG,h^F\right\>_{\H^F}^2=\sup_{\|h^F\|_{\H^F}\le1}\big|D_{h^F}G\big|^2\\
			&=\sup_{\|\dot{h}\|_{L^2([0,T];\d x)}\le1}\big|D_{K_F\dot{h}}G\big|^2=\sup_{\|\dot{h}\|_{L^2([0,T];\d x)} \le 1}\left\<DG,K_F\dot{h}\right\>_{L^2([0,T];\d x)}^2\\
			&\le\|DG\|_{L^2([0,T];\d x)}^2\sup_{\|\dot{h}\|_{L^2([0,T];\d x)}\le1}\|K_F\dot{h}\|_{L^2([0,T];\d x)}^2.
		\end{aligned}
	\end{equation}
	Now we need to estimate $\sup_{\|\dot{h}\|_{L^2([0,T];\d x)}\le1}\|K_F\dot{h}\|_{L^2([0,T];\d x)}^2$. In fact,
$$\begin{aligned}
\|K_F\dot{h}\|_{L^2([0,T];\d x)}^2&=\int_0^T\left(\int_0^tF(t,s)\dot{h}(s)\d s\right)^2\d t\\
&\le\int_0^T\left(\int_0^TF(t,s)^2\d s\int_0^T\dot{h}(s)^2\d s\right)\d t\\
&=\left(\int_0^T\int_0^TF(t,s)^2\d s\d t\right)\left(\int_0^T\dot{h}(s)^2\d s\right)\\
&=\|F\|_{L^2([0,T]^2;\d x\times \d x)}^2\|\dot{h}\|_{L^2([0,T];\d x)}^2.
\end{aligned}
$$
This implies conclusion holds.
\end{proof}

\begin{remark}
	In the case of fractional Brownian motion, Fan \cite{F15} obtained logarithmic Sobolev inequalities for the $L^2$-Dirichlet form. For $F$-Gaussian process with stronger condition, \cite{CSW} also obtained the Logarithmic Sobolev inequality for the three types of Dirichlet forms.
\end{remark}

\section{Bismut-Elworthy-Li's formula}
In this section, based on the quasi-invariance theorem \ref{thmQuasi}, we will establish Bismut-Elworthy-Li's formula with respect to Gaussian Wiener measure.

For each $x\in \R$, define
$$B_t^{F,x}=B_t^F+x.
$$
Let $\mu$ be the distribution of stochastic processes $B_\cdot^{F,x}$.

\begin{thm}[\textbf{Bismut-Elworthy-Li's formula}]
\label{thmBEL}
Assume that
$$V_T:=\int_0^T F(T,s)^2\d s>0.
$$
Then, for any function $g\in L^2(\R;\d x)$,
	\begin{equation}\notag
		\frac{\d}{\d x}\EE\left(g\left(B_T^{F,x}\right)\right)
		=\frac{1}{V_T}\EE\left(g\left(B_T^{F,x}\right)B_T^F\right).
	\end{equation}
\end{thm}

\begin{proof}
	For each $\varepsilon\in\R$ and $t\geq0$, define
	$$h^{F,\varepsilon}_t:=-K_F\left(\frac{\varepsilon F(T,\cdot)}{V_T}\right)(t)
	=-\frac{\varepsilon}{V_T}\int_0^tF(t,s)F(T,s)\d s.
	$$
	Then $h^{F,\varepsilon}\in \H^F$. Let $u^{F,\varepsilon}:=-\varepsilon F(T,\cdot)/V_T$ so that $K_Fu^{F,\varepsilon}=h^{F,\varepsilon}$. Denote
	$$B_t^{F,x,\varepsilon}=B_t^{F,x}+h_t^{F,\varepsilon},\quad t\geq0.
	$$
	Let $\P^{\varepsilon}$ be the distribution of stochastic processes $B_\cdot^{F,x,\varepsilon}$. According to Theorem \ref{thmQuasi}, we have
	\begin{equation}\notag
		\frac{\d\P^{\varepsilon}}{\d\P}=\beta_{-h^{F,\varepsilon}}=\exp\left\{-\int_0^T u^{F,\varepsilon}(s)\d B_s-\frac{\varepsilon^2}{2V_T}\right\}.
	\end{equation}
	By Theorem \ref{thmQuasi} again, the law of  $B_{[0,T]}^{F,x+\varepsilon,\varepsilon}$ under the probability measure $\P^{\varepsilon}$ is the same as $B_{[0,T]}^{F,x+\varepsilon}$ under $\P$, i.e. for any bounded $G$,
	\begin{equation}\notag
		\int_\Omega G\left(B_{[0,T]}^{F,x+\varepsilon,\varepsilon}\right)\frac{\d\P^{\varepsilon}}{\d\P}\d \P=\int_\Omega G\left(B_{[0,T]}^{F,x+\varepsilon}\right)\d \P.
	\end{equation}
	Since
	$$\begin{aligned} B_T^{F,x+\varepsilon,\varepsilon}&=B_T^{F,x+\varepsilon}+h_T^{F,\varepsilon}=x+\varepsilon+B_T^{F}+h_T^{F,\varepsilon}\\
	&=x+\varepsilon+B_T^{F}-\frac{\varepsilon}{V_T}\int_0^T F(T,s)^2\d s\\
	&=x+\varepsilon+B_T^{F}-\varepsilon\\
	&=B_T^{F,x},\end{aligned}
	$$
	in particular, taking $G(\gamma)=g(\gamma_T)$, we obtain that
	$$\int_\Omega g\left(B_T^{F,x}\right)\frac{\d\P^{\varepsilon}}{\d\P}\d \P=\int_\Omega g\left(B_T^{F,x+\varepsilon,\varepsilon}\right)\frac{\d\P^{\varepsilon}}{\d\P}\d \P=\int_\Omega g\left(B_T^{F,x+\varepsilon}\right)\d \P.
	$$

	Since bounded functions $g$ are dense in $L^2(\R;\d x)$, $\beta_{-h^{F,\varepsilon}}={\d\P^{\varepsilon}}/{\d\P} \in L^2(\P)$ and $g(B_T^{F,x}) \in L^2(\P)$ for any $x\in \R$ and any $g\in L^2(\R;\d x)$.
	Thus, for any $g\in L^2(\R;\d x)$
	\begin{equation}\notag
		\int_{\Omega}g\left(B_T^{F,x+\varepsilon}\right)\d\P=\int_{\Omega}g\left(B_T^{F,x}\right)\exp\left\{-\int_0^T u^{F,\varepsilon}(s)\d B_s-\frac{\varepsilon^2}{2V_T}\right\}\d\P.
	\end{equation}
	Taking the derivative with respect to $\varepsilon$ on both sides of the above equation, and letting $\varepsilon=0$, we immediately get
	\begin{equation}\notag
		\frac{\d}{\d x}\EE\left(g(B_T^{F,x})\right)
		=\frac{1}{V_T}\EE\left(g(B_T^{F,x})B_T^F\right).
	\end{equation}
\end{proof}
\begin{remark}
	Fan \cite[Theorem 3.1]{F} had proved a Bismut-Elworthy-Li's formula for SDEs driven by fractional Brownian motion $F=K_H$ with $H\geq 1/2$. Chen-Sun-Wu \cite[Theorem 6.1]{CSW} also proved another form of Bismut-Elworthy-Li's formula with $H\geq 1/2$.
\end{remark}

\end{document}